\documentclass[12pt,reqno]{amsart}
\usepackage{amsmath}
\usepackage{amssymb}
\usepackage{amsthm}

\newtheorem{theorem}{Theorem}[section]
\newtheorem{prop}[theorem]{Proposition}
\newtheorem{lem}[theorem]{Lemma}
\newtheorem*{cor}{Corollary}
\theoremstyle{definition}
\newtheorem{defn}[theorem]{Definition}
\theoremstyle{remark}
\newtheorem*{rem}{Remark}
\numberwithin{equation}{section}

\newcommand{\Gg}{\mathfrak{g}}    %Gothic
\newcommand{\Gh}{\mathfrak{h}}

\begin{document}

\title[Quantum Heisenberg group]
{Construction of a quantum Heisenberg group}

\author{Byung-Jay Kahng}
\date{}
\address{Department of Mathematics\\ University of Kansas\\
Lawrence, KS 66045}
\email{bjkahng@math.ku.edu}
%\subjclass[2000]{46L65, 46L51, 81R50}
%\keywords{}

\begin{abstract}
In this paper, we give a construction of a ($C^*$-algebraic) quantum
Heisenberg group.  This is done by viewing it as the dual quantum group
of the specific non-compact quantum group $(A,\Delta)$ constructed earlier
by the author.  Our definition of the quantum Heisenberg group is different
from the one considered earlier by Van Daele.  To establish our object of
study as a locally compact quantum group, we also give a discussion on its
Haar weight, which is no longer a trace.  In the latter part of the paper,
we give some additional discussion on the duality mentioned above.
\end{abstract}
\maketitle

{\sc Introduction.}
Among the simplest while useful of non-compact groups is the Heisenberg
Lie group $H$, which is two-step nilpotent.  Our goal in this paper is
to construct a version of a quantum Heisenberg group (i.\,e. a ``quantized
$C_0(H)$''), which would be an example of a non-compact, $C^*$-algebraic
quantum group.

 Our example is certainly not genuinely new.  Already in the early 90's,
Van Daele in \cite{VD} gave a construction of a quantum Heisenberg group,
which was actually one of the first examples of a non-compact quantum group
in the $C^*$-algebra setting.  Similar example but with a different approach
(via geometric quantization) was given by Szymczak and Zakrzewski \cite{SZ}.
Meanwhile, a dual counterpart to these examples was given by Rieffel \cite
{Rf5}.  Ours is different from these, but it is true that we were strongly
motivated by these early examples.

 In \cite{BJKp2}, we constructed a specific non-compact quantum group
$(A,\Delta)$, by deformation quantization of a certain non-linear Poisson
structure.  The construction was based on a generalization of Rieffel's
approach (as given in \cite{Rf3} and \cite {Rf5}).  And we saw in our
previous papers that $(A,\Delta)$ can be considered as a ``quantum
Heisenberg group algebra'' (i.\,e. a ``quantized $C^*(H)$'').  Naturally,
we are interested in its dual counterpart.  The dual quantum group, to be
denoted by $(\hat{A},\hat{\Delta})$ in the below, will be our main object
of study in this paper.  It will be our candidate to be a quantum Heisenberg
group.

 By general theory (for instance, see \cite{KuVa}), the dual object of
a locally compact quantum group is again a locally compact quantum group.
This means that the proof of $(\hat{A},\hat{\Delta})$ being a locally compact
quantum group is more or less automatic from the proof of $(A,\Delta)$ being
one.  For this reason, we did not find pressing needs for giving a separate
presentation on $(\hat{A},\hat{\Delta})$ until now, and we instead have been
only giving indications of its existence on several occasions in our previous
papers \cite{BJKp2}, \cite{BJKhj}, \cite{BJKppha}.  However, as we are trying
to develop some applications of these quantum groups (some of the program
were already carried out in \cite{BJKhj} and \cite{BJKppdress}), and also
when we try to construct the ``quantum double'' (work in preparation), it
became necessary to clarify the notion of our quantum Heisenberg group.

 It is true that ours is not one of the attention-grabbing examples.
But it is modestly interesting on its own, just as an ordinary Heisenberg
group is an interesting object of study in various branches of mathematics.
So in this article, we plan to carry out a careful construction of the
quantum Heisenberg group $(\hat{A},\hat{\Delta})$, including its non-tracial
Haar weight.  We will try to make the discussion as detailed as
possible, even if we may have to repeat some of our earlier results.
On the other hand, note that even for the case of the (simpler) quantum
Heisenberg group of Van Daele, so far no explicit discussion in the
$C^*$-algebra setting on the Haar weight has been given.  

 Here is a quick summary of how this paper is organized.  In section 1,
we review our example $(A,\Delta)$.  Although the information at the
Poisson--Lie group level played a significant role in its construction,
that angle will be de-emphasized here for the purpose of brevity.  Among
the useful tools that appear is the multiplicative unitary operator $U_A$.
The description of our quantum Heisenberg group $(\hat{A},\hat{\Delta})$
is given in section 2.  After giving a realization of the underlying
$C^*$-algebra $\hat{A}$, we will construct its quantum group structures,
including comultiplication, antipode, and Haar weight.  We will make a
point that $(\hat{A},\hat{\Delta})$ is reasonable to be considered as
a ``quantized $C_0(H)$''.

 In section 3, we give a light discussion on the duality between $(A,
\Delta)$ and $(\hat{A},\hat{\Delta})$.  And towards the end, we mention
some other related quantum Heisenberg group algebras and quantum Heisenberg
groups, namely the ``opposite'' and ``co-opposite'' versions of $(A,\Delta)$
and $(\hat{A},\hat{\Delta})$.  Most of the results here are straight from
the general theory, but several of these will be useful in our future
applications, including the quantum double construction.

\section{The Hilbert space ${\mathcal H}$. The quantum group $(A,\Delta)$.}

 Let $H$ be the $(2n+1)$-dimensional Heisenberg Lie group.  The underlying
space for this Lie group is $\mathbb{R}^{2n+1}$, and the multiplication
on it is defined by
$$
(x,y,z)(x',y',z')=\bigl(x+x',y+y', z+z'+\beta(x,y')\bigr),
$$
for $x,y,x',y'\in\mathbb{R}^n$ and $z,z'\in\mathbb{R}$.  Here $\beta(\ ,\ )$
is the usual inner product on $\mathbb{R}^n$, used here for a possible
future generalization.

 As indicated in the Introduction, we wish to obtain our {\em quantum Heisenberg
group\/} (``quantized $C_0(H)$''), as the dual object to the non-compact quantum
group $(A,\Delta)$ constructed earlier by the author.  So let us give here
a brief review of definitions involving $(A,\Delta)$.  See \cite{BJKp2},\cite
{BJKppha}, for more detailed discussion.

 To begin with, we need to establish our underlying Hilbert space 
${\mathcal H}$.  To do this, we first consider the Heisenberg group $H$
as a finite-dimensional vector space.  As above, typical elements in $H$
are written as $(x,y,z)$.  Next, let $H^*$ be the dual vector space of $H$,
whose typical elements will be written as $(p,q,r)$.  Note that in \cite
{BJKp2}, we used the notation $\Gg$ for the space $H^*$ and considered it
as $\Gg={\Gh}^*$, where $\Gh$ is the Lie algebra (so a vector space)
corresponding to $H$.  Since $H=\Gh$ as vector spaces (by virtue of being
nilpotent), this is equivalent.

 We take the natural Lebesgue measure $dxdydz$ on $H$, which would be the
Haar measure for the group $H$.  Whereas on $H^*$, we consider the dual
Plancherel Lebesgue measure $dpdqdr$, corresponding to the chosen Haar
measure on $H$.  Then we can define the Fourier transform from $L^2(H)$
to $L^2(H^*)$, as follows:
$$
({\mathcal F}\xi)(p,q,r)=\int_H\bar{e}(p\cdot x+q\cdot y+r\cdot z)
\xi(x,y,z)\,dxdydz.
$$
Here $\cdot$ denotes the dual pairing, and $e(\ )$ is the function defined
by $e(t)=e^{2\pi it}$.  So $\bar{e}(t)=e^{-2\pi it}$.  By our choice of
measures, the Fourier transform is a unitary operator whose inverse is the
following:
$$
({\mathcal F}^{-1}\zeta)(x,y,z)=\int_{H^*} e(p\cdot x+q\cdot y+r\cdot z)
\zeta(x,y,z)\,dpdqdr.
$$
The Fourier inversion theorem (the unitarity of the Fourier transform)
holds such that we have: ${\mathcal F}^{-1}({\mathcal F}\xi)=\xi$ and
${\mathcal F}({\mathcal F}^{-1}\zeta)=\zeta$, at the level of $L^2$-functions
as well as at the level of Schwartz functions.

 By Fourier transform, we can regard $L^2(H)$ and $L^2(H^*)$ as more
or less the same.  Actually, it is more convenient to work with the
$L^2$-functions in the $(x,y,r)$ variables, which we denote by
${\mathcal H}$.  That is, ${\mathcal H}=L^2(H/Z\times H^*/{Z^{\perp}})$,
where $Z=\bigl\{{\text {$(0,0,z)$'s}}\bigr\}$ in $H$.  By using the
partial Fourier transform in the third variable (defined similarly
as above), we can see that ${\mathcal H}$ is isomorphic to $L^2(H)$
(as well as to $L^2(H^*)$).  All our constructions will be carried out
over the Hilbert space ${\mathcal H}$.

 As a $C^*$-algebra, $A$ is isomorphic to the twisted crossed product
$C^*$-algebra $C^*\bigl(H/Z,C_0(H^*/{Z^{\perp}}),\sigma\bigr)$, with
the twisting given by a certain cocycle term $\sigma$.  To be more
precise, consider $\mathcal A$, which is the space of Schwartz functions
in the $(x,y,r)$ variables having compact support in the $r$ variable.
Clearly, ${\mathcal A}\subseteq C_0(H/Z\times H^*/{Z^{\perp}})$ as
well as ${\mathcal A}\subseteq{\mathcal H}$.  For $f,g\in{\mathcal A}$,
define:
$$
(L_f g)(x,y,r):=\int f(\tilde{x},\tilde{y},r)g(x-\tilde{x},y-\tilde{y},r)
\bar{e}\bigl[\eta_{\lambda}(r)\beta(\tilde{x},y-\tilde{y})\bigr]\,d\tilde{x}
d\tilde{y}.
$$

\begin{rem}
In the definition above,  $\lambda\in\mathbb{R}$ is a fixed constant, which
determines a certain non-linear Poisson structure when $\lambda\ne0$.  The
expression $\eta_{\lambda}(r)$ is defined such that $\eta_{\lambda}(r)=
\frac{e^{2\lambda r}-1}{2\lambda}$, which reflects the non-linear flavor.
When $\lambda=0$, we take $\eta_{\lambda=0}(r)=r$.  We are not planning
to explicitly mention the Poisson structure here.  But in section 1 of
\cite{BJKp2}, we gave a discussion on how it is related with a so-called
``classical $r$-matrix'' element.  Finally, the expression $\bar{e}\bigl
[\eta_{\lambda}(r)\beta(\tilde{x},y-\tilde{y})\bigr]$ is the cocycle
term, indicated by $\sigma$ above.
\end{rem}

 In this way, we define the ``regular representation'' $L$, and obtain the
$C^*$-algebra $A$ as the norm closure in ${\mathcal B}({\mathcal H})$ of
$L({\mathcal A})$.  We will, in many occasions, regard $f\in{\mathcal A}$
and $L_f$ as the same, and consider ${\mathcal A}$ as a (dense) subalgebra
of $A$.  Actually, ${\mathcal A}$ is a ${}^*$-subalgebra of $A$, whose
multiplication is given by $L_{f\times g}=L_f L_g$.  The involution
$f\mapsto f^*$ can be described by $L_{f^*}=(L_f)^*$.  We have the
following:
\begin{align}\label{(Aops)}
(f\times_A g)(x,y,r)&=\int f(\tilde{x},\tilde{y},r)g(x-\tilde{x},y-\tilde{y},r)
\bar{e}\bigl[\eta_{\lambda}(r)\beta(\tilde{x},y-\tilde{y})\bigr]\,d\tilde{x}
d\tilde{y}.  \notag \\
f^*(x,y,r)&=\bar{e}\bigl[\eta_{\lambda}(r)\beta(x,y)\bigr]\overline{f(-x,-y,r)}.
\end{align}
See Propositions 2.8 and 2.9 of \cite{BJKp2}.

 In \cite{BJKppha} (in Proposition 2.2), we gave another characterization
of $A$, in terms of a certain unitary operator $U_A\in{\mathcal B}
({\mathcal H}\otimes{\mathcal H})$.  Namely,
\begin{equation}\label{(Lomega)}
A=\overline{\bigl\{(\omega\otimes\operatorname{id})(U_A):\omega\in
{\mathcal B}({\mathcal H})_*\bigr\}}^{\|\ \|},
\end{equation}
where the $L(\omega)=(\omega\otimes\operatorname{id})(U_A)\in{\mathcal B}
({\mathcal H})$ are the ``left slices'' of $U_A$ by the linear forms
$\omega\in{\mathcal B}({\mathcal H})_*$.  The unitary operator $U_A$
is ``multiplicative'' (in the sense of Baaj and Skandalis \cite{BS}),
and is defined by
\begin{align}
U_A\xi(x,y,r,x',y',r')&=(e^{-\lambda r'})^n\bar{e}\bigl[\eta_{\lambda}
(r')\beta(e^{-\lambda r'}x,y'-e^{-\lambda r'}y)\bigr] \notag \\
&\quad\xi(e^{-\lambda r'}x,e^{-\lambda r'}y,r+r',
x'-e^{-\lambda r'}x,y'-e^{-\lambda r'}y,r').  \notag
\end{align}

 Using the multiplicative unitary operator $U_A$, we can also define
the {\em comultiplication\/} $\Delta:A\to M(A\otimes A)$.  For $a\in A$,
we have:
$$
\Delta a=U_A(a\otimes1){U_A}^*.
$$
The comultiplication is a non-degenerate ${}^*$-homomorphism satisfying
the coassociativity condition: $(\Delta\otimes\operatorname{id})(\Delta a)
=(\operatorname{id}\otimes\Delta)(\Delta a)$.  In case $f\in{\mathcal A}$,
the equation $\Delta(L_f)=(L\otimes L)_{\Delta f}$ gives us the following:
\begin{align}\label{(Acoproduct)}
&{\Delta}f(x,y,r;x',y',r')  \notag \\
&=\int f(x',y',r+r')\bar{e}\bigl[\tilde{p}\cdot(e^{\lambda r'}x'-x)
+\tilde{q}\cdot(e^{\lambda r'}y'-y)\bigr]\,d\tilde{p}d\tilde{q},
\end{align}
which is a Schwartz function having compact support in $r$ and $r'$.

 We have further shown in \cite{BJKp2} and \cite{BJKppha} that the
$C^*$-bialgebra $(A,\Delta)$, together with the additional structures
on it, namely counit $\varepsilon$, antipode $S$, and the Haar weight
$\varphi_A$, becomes a locally compact ($C^*$-algebraic) quantum group,
in the sense of Kustermans and Vaes \cite{KuVa}.  In particular, a
rigorous discussion on the definition of $\varphi_A$ as a $C^*$-algebra
weight and its left invariance property was given in section 3 of
\cite{BJKppha}.

\begin{rem}
We have been arguing in our earlier papers that $(A,\Delta)$ is, in a
sense, a ``quantum Heisenberg group algebra'' (For instance, see \cite
{BJKhj}, where we studied its representation theory.).  To make a brief
case here, let $\lambda=0$ (so $\eta_{\lambda}(r)=r$).  Let us now
re-write the expression for $L_f g$ at the level of $(x,y,z)$ variables,
by using the partial Fourier transform in the third variable and by
using the Fourier inversion theorem.  Then we have (for convenience,
we are not distinguishing a function from its partial Fourier transform):
\begin{align}
&(L_f g)(x,y,z)  \notag \\
&=\int f(\tilde{x},\tilde{y},\tilde{z})\bar{e}(\tilde{z}r)g(x-\tilde{x},
y-\tilde{y},\hat{z})\bar{e}(\hat{z}r)\bar{e}\bigl[r\beta(\tilde{x},
y-\tilde{y})\bigr]e(rz)\,d\tilde{z}d\hat{z}d\tilde{x}d\tilde{y}dr  \notag \\
&=\int f(\tilde{x},\tilde{y},\tilde{z})g\bigl(x-\tilde{x},y-\tilde{y},
z-\tilde{z}-\beta(\tilde{x},y-\tilde{y})\bigr)\,d\tilde{x}d\tilde{y}d\tilde{z}.
\notag
\end{align}
This is just the left regular representation of $C^*(H)$, given by the
convolution product.  The involution can be also realized as the natural
one on the convolution algebra.
\end{rem}

\section{The quantum Heisenberg group $(\hat{A},\hat{\Delta})$.}

 Since $(A,\Delta)$ can be viewed as a ``quantized $C^*(H)$'', it is
natural to consider its dual object as the candidate for the quantum
Heisenberg group.  Suggested by the general theory of locally compact
quantum groups \cite{KuVa}, and taking advantage of the theory of
multiplicative unitary operators \cite{BS}, we will define $(\hat{A},
\hat{\Delta})$ in terms of our fundamental multiplicative unitary
operator $U_A$.

\begin{defn}\label{Ahat}
Consider the ``right slices'' of $U_A$, which are the operators $\rho(\omega)
\in{\mathcal B}({\mathcal H})$ defined by $\rho(\omega)=(\operatorname{id}
\otimes\omega)(U_A)$, for $\omega\in{\mathcal B}({\mathcal H})_*$.  We will
define $\hat{A}$ as the $C^*$-algebra generated by the $\rho(\omega)$:
$$
\hat{A}=\overline{\bigl\{(\operatorname{id}\otimes\omega)(U_A):\omega\in
{\mathcal B}({\mathcal H})_*\bigr\}}^{\|\ \|}.
$$
For a typical element $b\in\hat{A}$, define $\hat{\Delta}b$ by $\hat
{\Delta}b={U_A}^*(1\otimes b)U_A$.  In this way, we obtain the {\em 
comultiplication\/} $\hat{\Delta}:\hat{A}\to M(\hat{A}\otimes\hat{A})$,
which is a non-degenerate $C^*$-homomorphism satisfying the coassociativity.
\end{defn}

 The general theory of multiplicative unitary operators assures us that
$(\hat{A},\hat{\Delta})$ is a $C^*$-bialgebra (or a quantum semigroup)
\cite{BS}, \cite{Wr7}, \cite{KuVa}.  But to be more specific and to be
more accessible in our future applications, let us give here an explicit
realization of the $C^*$-algebra $\hat{A}$:

\begin{prop}\label{rho}
Let $\hat{\mathcal A}$ be the space of Schwartz functions in the $(x,y,r)$
variables having compact support in the $r$ variable.  For $f\in\hat
{\mathcal A}$, define the operator ${\rho}_f\in{\mathcal B}({\mathcal H})$
by
$$
({\rho}_f\zeta)(x,y,r)=\int(e^{\lambda\tilde{r}})^n f(x,y,\tilde{r})
\zeta(e^{\lambda\tilde{r}}x,e^{\lambda\tilde{r}}y,r-\tilde{r})\,d\tilde{r}.
$$
Then the $C^*$-algebra $\hat{A}$ is generated by the operators ${\rho}_f$.
\end{prop}

\begin{proof}
Let us work with the standard notation $\omega_{\xi,\eta}$, where $\xi,\eta
\in{\mathcal H}$.  It is defined by $\omega_{\xi,\eta}(T)=\langle T\xi,\eta
\rangle$, for $T\in{\mathcal B}({\mathcal H})$.  It is well known that linear
combinations of the $\omega_{\xi,\eta}$ are (norm) dense in ${\mathcal B}
({\mathcal H})_*$.  So consider $(\operatorname{id}\otimes\omega_{\xi,\eta})
(U_A)\in{\mathcal B}({\mathcal H})$. Without loss of generality, we can assume
that $\xi$ and $\eta$ are continuous functions having compact support.  Then
for $\zeta\in{\mathcal H}$, we have:
\begin{align}
\bigl((\operatorname{id}\otimes\omega_{\xi,\eta})(U_A)\bigr)\zeta(x,y,r)
&=\int\bigl(U_A(\zeta\otimes\xi)\bigr)(x,y,r;\tilde{x},\tilde{y},\tilde{r})
\overline{\eta(\tilde{x},\tilde{y},\tilde{r})}\,d\tilde{x}d\tilde{y}d\tilde{r}
\notag \\
&=\int (e^{\lambda\tilde{r}})^n f(x,y,\tilde{r})\zeta(e^{\lambda\tilde{r}}x,
e^{\lambda\tilde{r}}y,r-\tilde{r})\,d\tilde{r},  \notag
\end{align}
where
$$
f(x,y,\tilde{r})=\int\bar{e}\bigl[\eta_{\lambda}(\tilde{r})\beta
(x,y-e^{-\lambda\tilde{r}}\tilde{y})\bigr]\xi(\tilde{x}-e^{\lambda
\tilde{r}}x,\tilde{y}-e^{\lambda\tilde{r}}y,-\tilde{r})\overline
{\eta(\tilde{x},\tilde{y},-\tilde{r})}\,d\tilde{x}d\tilde{y}.
$$
Since $\xi$ and $\eta$ are $L^2$-functions, $f$ is a continuous function.
It would also have compact support.  Meanwhile, since the choice of $\xi$
and $\eta$ is arbitrary and since the $\omega_{\xi,\eta}$ are dense in
${\mathcal B}({\mathcal H})_*$, we can see that the collection of the
$f$ will form a total set in the space of continuous functions in the
$(x,y,r)$ variables having compact support.  It follows that we have:
$$
\overline{\rho(\hat{\mathcal A})}^{\|\ \|}=
\overline{\bigl\{(\operatorname{id}\otimes\omega)(U_A):\omega\in
{\mathcal B}({\mathcal H})_*\bigr\}}^{\|\ \|}=\hat{A}.
$$
\end{proof}

 As in the case of $A$, we will often regard the functions $f\in
\hat{\mathcal A}$ as the same as the operators $\rho_f\in\hat{A}$.
In this way, $\hat{\mathcal A}$ is considered as a ${}^*$-subalgebra
of $\hat{A}$.  The multiplication on it is defined by $\rho_{f\times g}
=\rho_f\rho_g$:
\begin{equation}\label{(multiplication)}
(f\times_{\hat{A}} g)(x,y,r)=\int f(x,y,\tilde{r})g(e^{\lambda\tilde{r}}x,
e^{\lambda\tilde{r}}y,r-\tilde{r})\,d\tilde{r}.
\end{equation}
While, the involution on $\hat{\mathcal A}$ is given by $\rho_{f^*}
=(\rho_f)^*$:
\begin{equation}\label{(involution)}
f^*(x,y,r)=\overline{f(e^{\lambda r}x,e^{\lambda r}y,-r)}.
\end{equation}

\begin{rem}
We can see that as a $C^*$-algebra, $\hat{A}\cong C_0(\mathbb{R}^{2n})
\rtimes_{\alpha}\mathbb{R}$, which is a crossed product algebra (together
with the action of $\mathbb{R}$ onto $\mathbb{R}^{2n}$ given by $\alpha(r):
(x,y)\mapsto(e^{\lambda r}x,e^{\lambda r}y)$).  Although our construction
of $\hat{A}$ here is given indirectly by relying on the duality, certainly
there is a direct way of constructing $(\hat{A},\hat{\Delta})$, giving a
Hopf $C^*$-algebra structure on a crossed product algebra (It is actually
easier than the case of $(A,\Delta)$ in \cite{BJKp2}.  See also the more
general approach described in \cite{VV}.).  Meanwhile, we also note that
when the parameter $\lambda=0$, we have: $\hat{A}_{\lambda=0}\cong C_0
(\mathbb{R}^{2n+1})=C_0(H)$, giving us a mild justification that $\hat{A}$
is a good candidate to become a ``quantized $C_0(H)$''.
\end{rem}

 Turning our attention to the coalgebra structure, let us give here
the description of the comultiplication $\hat{\Delta}$, at the level
of functions:

\begin{prop}\label{comultiplication}
For $f\in\hat{\mathcal A}$, let $\hat{\Delta}f$ be the Schwartz function
in the $(x,y,r;x',y',r')$ variables having compact support in $r$ and $r'$,
given by
\begin{align}
&\hat{\Delta}f(x,y,r;x',y',r')   \notag \\
&=\int f(x+x',y+y',\tilde{r})e\bigl[\eta_{\lambda}(\tilde{r})\beta(x,y')
\bigr]e\bigl[\tilde{r}(z+z')\bigr]\bar{e}[zr+z'r']\,d\tilde{r}dzdz'.
\notag
\end{align}
The map $f\mapsto\hat{\Delta}f$ extends to the map $\hat{\Delta}:
\hat{A}\to M(\hat{A}\otimes\hat{A})$, which is the comultiplication
given in Definition \ref{Ahat}.
\end{prop}

\begin{proof}
For $f\in\hat{\mathcal A}$ and for $\xi\in{\mathcal H}\otimes{\mathcal H}$,
we have:
\begin{align}
&{U_A}^*(1\otimes\rho_f)U_A\xi(x,y,r;x',y',r')  \notag \\
&=\int (e^{2\lambda\tilde{r}})^n f(x+x',y+y',\tilde{r})
e\bigl[\eta_{\lambda}(\tilde{r})\beta(x,y')\bigr]
\xi(e^{\lambda\tilde{r}}x,e^{\lambda\tilde{r}}y,r-\tilde{r};
e^{\lambda\tilde{r}}x',e^{\lambda\tilde{r}}y',r'-\tilde{r})\,d\tilde{r}.
\notag
\end{align}
By straightforward calculation, we can check without difficulty that
$(\rho\otimes\rho)_{\hat{\Delta}f}\xi(x,y,r;x',y',r')={U_A}^*
(1\otimes\rho_f)U_A\xi(x,y,r;x',y',r')$.  This means that $f\mapsto\hat
{\Delta}f$ coincides with the comultiplication $\rho_f\mapsto\hat{\Delta}
(\rho_f)$, at the level of the dense subalgebra.  It obviously extends to
the comultiplication on $\hat{A}$.
\end{proof}

 We will skip the proofs of the various properties of $\hat{\Delta}$,
referring instead to general theory.  For instance, the coassociativity
of $\hat{\Delta}$ follows from the unitary operator $U_A$ being
multiplicative.  Meanwhile, we see again that when the parameter
$\lambda=0$, the expression at the level of functions in the $(x,y,z;
x',y',z')$ variables for $\hat{\Delta}f$, obtained by using partial
Fourier transform, is just:
$$
\hat{\Delta}f(x,y,z;x',y',z')=f\bigl(x+x',y+y',z+z'+\beta(x,y')\bigr),
$$
recovering the usual comultiplication on $C_0(H)$.

 Next, let us consider the {\em antipode\/} $\hat{S}$.  The main result
is summarized below:

\begin{prop}\label{antipode}
For $f\in\hat{\mathcal A}$, let $\hat{S}(f)$ be the function in $\hat
{\mathcal A}$ defined by
$$
\bigl(\hat{S}(f)\bigr)(x,y,r)=\bar{e}\bigl[\eta_{\lambda}(r)\beta(x,y)\bigr]
f(-e^{\lambda r}x,-e^{\lambda r}y,-r).
$$
Then $\hat{S}$ can be extended to the anti-automorphism $\hat{S}:\hat{A}
\to\hat{A}$.  It is the antipodal map, satisfying: $\hat{S}\bigl(\hat{S}
(b)^*\bigr)^*=b$ and $(\hat{S}\otimes\hat{S})(\hat{\Delta}b)=\chi\bigl(
\hat{\Delta}(\hat{S}(b))\bigr)$, where $\chi$ denotes the flip.  We
also have: $\hat{S}^2\equiv\operatorname{Id}$.
\end{prop}

\begin{proof}
The definition is suggested by \cite{BS} and \cite{Wr7}.  It is equivalent
to the map $\hat{S}:(\operatorname{id}\otimes\omega)(U_A)\mapsto(\operatorname
{id}\otimes\omega)({U_A}^*)$, for $\omega\in{\mathcal B}({\mathcal H})_*$.
Using same kind of the technique we used in the proof of Proposition \ref{rho},
we could obtain the above expression for $\hat{S}(f)$.  As before, this
should be interpreted as $\hat{S}(\rho_f)=\rho_{\hat{S}(f)}$.  Meanwhile,
a straight calculation shows that $\hat{S}$ can be  equivalently written as
$\hat{S}(b)=Jb^*J$, for $b\in\hat{A}$, where $J$ is the anti-unitary operator
defined by
$$
J\xi(x,y,r)=\bar{e}\bigl[\eta_{\lambda}(r)\beta(x,y)\bigr]
\overline{\xi(-x,-y,r)}.
$$
Due to this characterization, the remaining properties are easy to verify.
\end{proof}

\begin{rem}
The notation for the operator $J$ introduced in the proof is motivated by
the modular theory, and it is essentially the involution on ${\mathcal A}$
(as defined in \eqref{(Aops)}).  Indeed, remembering that the space
${\mathcal A}$ is dense in ${\mathcal H}$ with respect to the Hilbert
space norm, we see that $J$ is just the extension of the map $f\mapsto f^*$
in ${\mathcal A}$.
\end{rem}

 The correct formulation of $\hat{S}$ being the legitimate antipode relies
on the existence of an appropriate Haar weight (to be constructed shortly).
However, we may still point out that if $\lambda=0$, the antipode at the
level of functions in the $(x,y,z)$ variables is just:
$$
\bigl(\hat{S}(f)\bigr)(x,y,z)=f\bigl(-x,-y,-z+\beta(x,y)\bigr)
=f\bigl((x,y,z)^{-1}\bigr).
$$
Here again, we used the partial Fourier transform.

 Since $\hat{S}$ is already an anti-automorphism such that $\hat{S}^2\equiv
\operatorname{Id}$, its ``polar decomposition'' (in the sense of \cite{MN}
and \cite{KuVa}) is trivial: That is, $\hat{S}=\hat{R}$ (the ``unitary
antipode''), and $\hat{\tau}\equiv\operatorname{Id}$ (the ``scaling group''). 
Meanwhile, $\hat{S}^2\equiv\operatorname{Id}$ suggests that our example
will be a kind of a ``Kac $C^*$-algebra'' \cite{Va}, \cite{ES}, which is
expected since $(A,\Delta)$ was one such.

 One remaining important structure to be constructed is the Haar weight.
We begin with the linear functional $\hat{\varphi}$ defined at the dense
function algebra level (i.\,e. on $\hat{\mathcal A}$), motivated by the
Lebesgue measure on $H$:

\begin{prop}\label{Haarfunctional}
On $\hat{\mathcal A}$, define the linear functional $\hat{\varphi}$ by
$$
\hat{\varphi}(f)=\int f(x,y,0)\,dxdy.
$$
Then $\hat{\varphi}$ defined as above is a faithful, positive linear
functional.  It is also unimodular, in the sense that $\hat{\varphi}
\circ\hat{S}=\hat{\varphi}$.
\end{prop}

\begin{proof}
Suppose $F\in\hat{\mathcal A}$ is a typical positive element such that $\rho_F
=(\rho_f)(\rho_f)^*$ for some $f\in\hat{\mathcal A}$.  Then we have:
\begin{align}
\hat{\varphi}(F)&=\hat{\varphi}(f\times_{\hat{A}}f^*)=\int(f\times_{\hat{A}}f^*)
(x,y,0)\,dxdy  \notag \\
&=\int f(x,y,\tilde{r})f^*(e^{\lambda\tilde{r}}x,e^{\lambda\tilde{r}}y,
0-\tilde{r})\,d\tilde{r}dxdy   \notag \\
&=\int f(x,y,\tilde{r})\overline{f(x,y,\tilde{r})}\,d\tilde{r}dxdy=\|f\|_2^2.
\notag
\end{align}
From this, the first part of the proposition is immediate.  Meanwhile, for
an arbitrary element $f\in\hat{\mathcal A}$, we have:
$$
\hat{\varphi}\bigl(\hat{S}(f)\bigr)=\int\bigl(\hat{S}(f)\bigr)(x,y,0)\,dxdy
=\int f(-x,-y,0)\,dxdy=\hat{\varphi}(f),
$$
giving us the proof of the unimodularity.
\end{proof}

 We need to find a $C^*$-algebra weight extending this linear functional.
The following steps are more or less the same ones we took in \cite{BJKppha}
(One difference is that $\hat{\varphi}$ is no longer a trace.).  First,
let us consider the GNS construction associated with $\hat{\varphi}$.  We
see below that the ``regular representation'' $\rho$ of $\hat{\mathcal A}$
we have been using is essentially the GNS representation:

\begin{prop}
Let $\Gamma:\hat{\mathcal A}\to{\mathcal H}$ be defined by
$\Gamma(f)(x,y,r):=(e^{\lambda r})^n f(x,y,r)$.  Then for $f,g\in\hat
{\mathcal A}$, we have:
$$
\bigl\langle\Gamma(f),\Gamma(g)\bigr\rangle_{\mathcal H}=\hat{\varphi}
(g^*\times_{\hat{A}}f),
$$
where $\langle\ ,\ \rangle_{\mathcal H}$ is the inner product on ${\mathcal H}$,
conjugate in the second place.  From this, we see that $\Gamma$ gives
the Hilbert space isomorphism between ${\mathcal H}_{\hat{\varphi}}$
and ${\mathcal H}$, where ${\mathcal H}_{\hat{\varphi}}$ is the GNS
Hilbert space for $\hat{\varphi}$.  Meanwhile, consider the non-degenerate
${}^*$-representation $\pi_{\hat{\varphi}}:\hat{\mathcal A}\to{\mathcal B}
({\mathcal H})$, given by $\bigl(\pi_{\hat{\varphi}}(f)\bigr)\bigl(\Gamma(g)
\bigr):=\Gamma(f\times_{\hat{A}}g)$.  It turns out that $\pi_{\hat{\varphi}}$
coincides with the representation $\rho$.
\end{prop}

\begin{proof}
For $f,g\in\hat{\mathcal A}$,
\begin{align}
\hat{\varphi}(g^*\times_{\hat{A}}f)&=\int g^*(x,y,\tilde{r})
f(e^{\lambda\tilde{r}}x,e^{\lambda\tilde{r}}y,0-\tilde{r})\,d\tilde{r}dxdy
\notag \\
&=\int\overline{g(e^{\lambda\tilde{r}}x,e^{\lambda\tilde{r}}y,-\tilde{r})}
f(e^{\lambda\tilde{r}}x,e^{\lambda\tilde{r}}y,-\tilde{r})\,d\tilde{r}dxdy
\notag \\
&=\int(e^{2\lambda\tilde{r}})^n\overline{g(x,y,\tilde{r})}
f(x,y,\tilde{r})\,d\tilde{r}dxdy=\bigl\langle\Gamma(f),\Gamma(g)\bigr
\rangle_{\mathcal H}.   \notag
\end{align}
Since the GNS Hilbert space ${\mathcal H}_{\hat{\varphi}}$ is obtained by
completing $\hat{\mathcal A}$ with respect to the inner product $(f,g)
\mapsto\hat{\varphi}(g^*\times_{\hat{A}}f)$, we see easily that $\Gamma$
(now extended to ${\mathcal H}_{\hat{\varphi}}$) provides the Hilbert
space isomorphism $\Gamma:{\mathcal H}_{\hat{\varphi}}\cong{\mathcal H}$.
The representation $\pi_{\hat{\varphi}}$ being non-degenerate is immediate,
remembering that $\Gamma(\hat{\mathcal A})$ is dense in ${\mathcal H}$.
Now to learn about the representation $\pi_{\hat{\varphi}}$, consider
$f,g\in\hat{\mathcal A}$.  Let us write $\zeta =\Gamma(g)\in{\mathcal H}$.
Then:
\begin{align}
\bigl(\pi_{\hat{\varphi}}(f)\bigr)\zeta(x,y,r)&=\bigl(\Gamma(f\times_{\hat{A}}
g)\bigr)(x,y,r)=(e^{\lambda r})^n(f\times_{\hat{A}}g)(x,y,r)  \notag \\
&=\int(e^{\lambda r})^n f(x,y,\tilde{r})g(e^{\lambda\tilde{r}}x,
e^{\lambda\tilde{r}}y,r-\tilde{r})\,d\tilde{r}  \notag \\
&=\int(e^{\lambda\tilde{r}})^n f(x,y,\tilde{r})\zeta(e^{\lambda\tilde{r}}x,
e^{\lambda\tilde{r}}y,r-\tilde{r})\,d\tilde{r}  \notag \\
&=(\rho_f\zeta)(x,y,r),  \notag
\end{align}
recovering the representation $\rho$.
\end{proof}

 By (essential) uniqueness of GNS construction, we see from the
above proposition that $({\mathcal H},\Gamma,\pi_{\hat{\varphi}}=\rho)$
is the GNS triple associated with $\hat{\varphi}$.  The consequence is
that the algebra $\hat{\mathcal A}$ (to be more precise, $\Gamma(\hat
{\mathcal A})\subseteq{\mathcal H}$) is a ``left Hilbert algebra''
(See literature on modular theory \cite{Cm2}, \cite{St}.).  One detail
to note is that the involution on $\hat{\mathcal A}$ is not isometric
with respect to the inner product (This reflects the fact that the
functional $\hat{\varphi}$ is not a trace.).  But it is still closable.

 We denote by $\hat{T}$ the closure of the involution on $\hat{\mathcal A}$.
Then it is a closed, anti-linear map on ${\mathcal H}$ having $\Gamma(\hat
{\mathcal A})$ as a core for $\hat{T}$, such that $\hat{T}\bigl(\Gamma(f)\bigr)
=\Gamma(f^*)$.  By a simple calculation, we have:
$$
\hat{T}\zeta(x,y,r)=(e^{2\lambda r})^n\overline{\zeta(e^{\lambda r}x,
e^{\lambda r}y,-r)}.
$$
We have the polar decomposition: $\hat{T}=\hat{J}{\hat{\nabla}}^{\frac{1}{2}}$,
where $\hat{\nabla}={\hat{T}}^*\hat{T}$ is the ``modular operator'' and $\hat{J}$
is an anti-unitary operator.  They are given as follows:
$$
\hat{\nabla}f(x,y,r)=(e^{-2\lambda r})^n f(x,y,r),\quad
\hat{J}f(x,y,r)=(e^{\lambda r})^n\overline{f(e^{\lambda r}x,
e^{\lambda r}y,-r)}.
$$
Note here that $\hat{J}$ is exactly the anti-unitary operator which
we used in our definition of the antipode $S$ for $(A,\Delta)$, given
in \cite{BJKppha}.  Compare this with the remark we made following
Proposition \ref{antipode}, pointing out the relationship between
the operator $J$ and the antipode $\hat{S}$ of $(\hat{A},\hat{\Delta})$.
This aspect is one of many useful relationships between the (mutually
dual) algebras ${\mathcal A}$ and $\hat{\mathcal A}$.  See \cite{MN}
and \cite {KuVa}.

 Since we have a left Hilbert algebra structure on $\hat{\mathcal A}$, we
can follow the standard modular theory (\cite{Cm2}, \cite{St})) to obtain
a $C^*$-algebra weight extending $\hat{\varphi}$.  The modular operator
$\hat{\nabla}$ plays an important role in the formulation of the KMS
property.

\begin{theorem}\label{c*weight}
There is a faithful, lower semi-continuous weight $\hat{\varphi}_{\hat{A}}$
on the $C^*$-algebra $\hat{A}$, extending the linear functional $\hat{\varphi}$.
It is also a KMS weight: With respect to the (norm-continuous) one-parameter
group of automorphisms $\hat{\sigma}$ given by $\hat{\sigma}_t(b)=\hat
{\nabla}^{it}b\hat{\nabla}^{-it}$, we have:
\begin{align}
&\hat{\varphi}_{\hat{A}}\circ\hat{\sigma}_t=\hat{\varphi}_{\hat{A}},\quad
{\text { for all $t\in\mathbb{R}$,}}  \notag \\
&\hat{\varphi}_{\hat{A}}(b^*b)=\hat{\varphi}_{\hat{A}}\bigl(\hat{\sigma}_{i/2}
(b){\hat{\sigma}_{i/2}(b)}^*\bigr),\quad {\text { for all $b\in D(\hat
{\sigma}_{i/2})$.}}   \notag
\end{align}
\end{theorem}

\begin{rem}
The notion of ``KMS weight'' we are using above is due to Kustermans \cite{Ku1},
which is actually equivalent to the original notion given by Combes \cite{Cm2}.
Since the weight $\hat{\varphi}_{\hat{A}}$ extends the functional $\hat{\varphi}$
on $\hat{\mathcal A}$, it is densely defined, giving us a ``proper'' KMS weight.
Refer the discussion in \S1.1 of \cite{BJKppha} or literature on weight theory
\cite{Cm1}, \cite{Cm2}, \cite{St}.
\end{rem}

\begin{proof}
The (non-degenerate) representation $\pi_{\hat{\varphi}}(=\rho)$ generates
the von Neumann algebra $M_{\hat{A}}=\rho(\hat{\mathcal A})''$ in ${\mathcal B}
({\mathcal H})$.  Since $\hat{\mathcal A}$ is a left Hilbert algebra, there
is a standard way of defining a faithful, semi-finite, normal weight on
$M_{\hat{A}}$ (See Theorem 2.11 of \cite{Cm2}.  See also the discussion
we made in Theorem 3.6 of \cite{BJKppha}.).  We then obtain our weight
$\hat{\varphi}_{\hat{A}}$, by restricting this normal weight to the
$C^*$-algebra $\hat{A}=\overline{\rho(\hat{\mathcal A})}^{\|\ \|}\bigl
(\subseteq\rho(\hat{\mathcal A})''=M_{\hat{A}}\bigr)$.  Because of the way
it is constructed, it is not difficult to see that $\hat{\varphi}_{\hat{A}}$
extends the functional $\hat{\varphi}$ and is faithful.

The lower semi-continuity and the KMS property of $\hat{\varphi}_{\hat{A}}$
is a consequence of the fact that it is obtained from a normal weight at the
von Neumann algebra level. In our case, the modular automorphism group is
such that $\hat{\mathcal A}$ forms a core for the $\hat{\sigma}_t$ and that
for $f\in{\mathcal A}$, we have:
$$
\bigl(\hat{\sigma}_t(f)\bigr)(x,y,r)=(e^{-2\lambda r it})^n f(x,y,r).
$$
As before, this is interpreted as $\rho_{\hat{\sigma}_t(f)}=\hat{\sigma}_t
(\rho_f)=\hat{\nabla}^{it}\rho_f\hat{\nabla}^{-it}$.  To verify the KMS
property, we just choose $f\in\hat{\mathcal A}$ can calculate.  Since
$\bigl(\hat{\sigma}_{i/2}(f)\bigr)(x,y,r)=(e^{\lambda r})^n f(x,y,r)$,
We have:
$$
\hat{\varphi}\bigl(\hat{\sigma}_{i/2}(f){\hat{\sigma}_{i/2}(f)}^*\bigr)
=\int(e^{2\lambda\tilde{r}})^n f(x,y,\tilde{r})\overline{f(x,y,\tilde{r})}
\,d\tilde{r}dxdy=\hat{\varphi}(f^*\times_{\hat{A}}f).
$$
Verification of $\hat{\varphi}\bigl(\hat{\sigma}_t(f)\bigr)=\hat{\varphi}(f)$
is also straightforward.
\end{proof}

 For $\hat{\varphi}_{\hat{A}}$ to be considered as the legitimate Haar weight
(as well as to complete the discussion that $(\hat{A},\hat{\Delta})$ is a
locally compact quantum group), we need to establish its (left) invariance
property.  This will be done following the idea suggested by Van Daele
\cite{VDaxb}, \cite{VDoamp} (See also our discussion in section 3 of \cite
{BJKppha}.).  We first begin with a lemma.

\begin{lem}
Let $M_A$ be the enveloping von Neumann algebra of $A$ (That is, $M_A=
L({\mathcal A})''$.), while $M_{\hat{A}}\bigl(=\rho(\hat{\mathcal A})''\bigr)$
is the enveloping von Neumann algebra of $\hat{A}$ as appeared in the proof
of Theorem \ref{c*weight}.  Then we have:
\begin{enumerate}
\item $U_A\in M_{\hat{A}}\otimes M_A\subseteq{\mathcal B}({\mathcal H}
\otimes{\mathcal H})$.
\item $M_A\cap M_{\hat{A}}=\mathbb{C}1$.
\item The linear space $M_A M_{\hat{A}}$ is $\sigma$-strongly dense in
${\mathcal B}({\mathcal H})$.
\end{enumerate}
\end{lem}

\begin{proof}
The first statement follows from general theory of multiplicative unitary
operators.  It is also true that we have: $U_A\in M(\hat{A}\otimes A)$.
For the next two statements, we may follow Proposition 2.5 of \cite{VDoamp}.
\end{proof}

 The main strategy suggested by Van Daele is that there exists a faithful,
semi-finite, normal weight $\nu$ on ${\mathcal B}({\mathcal H})$ such that
at least formally, $\nu(ab)={\varphi}(a)\hat{\varphi}(b)$, for $a\in M_A$,
$b\in M_{\hat{A}}$.  Note here that for convenience, we are using the notation
$\varphi$ and $\hat{\varphi}$ for the weights.  As long as there is not
going to be confusion between the linear functionals and the weights, we
will often use the simpler notation.

\begin{prop}\label{nutrace}
On ${\mathcal B}({\mathcal H})$, consider the linear functional $\nu:=
\operatorname{Tr}$.  Then $\nu$ is a faithful, semi-finite, normal weight on
${\mathcal B}({\mathcal H})$ such that for $a\in{\mathfrak N}_{{\varphi}}$
and $b\in{\mathfrak N}_{\hat{\varphi}}$,
$$
\nu(a^*b^*ba)=\varphi(a^*a)\hat{\varphi}(b^*b).
$$
\end{prop}

\begin{rem}
The notations ${\mathfrak N}_{{\varphi}}$ and ${\mathfrak N}_{\hat{\varphi}}$
are the standard ones used in weight theory, which just ensure that the
expression in the right hand side does make sense (i.\,e. finite).  This
result is actually quite general in nature, although $\nu$ should be in
general a certain ``weighted trace'' instead of being the regular trace
(See Definition 2.6 of \cite{VDoamp}.).  The reason why the regular trace
works in our case has to do with the fact that $\varphi$ is a tracial weight
on $M_A$ (as shown in section 3 of \cite{BJKppha}).
\end{rem}

\begin{proof}
Let us pick two elements at the dense function algebra level, namely
$a=L_a\in{\mathcal A}(\subseteq M_A)$ and $b=\rho_b\in\hat{\mathcal A}
(\subseteq M_{\hat{A}})$.  Then by definition of $L_a$ and $\rho_b$
given in earlier sections, we have:
\begin{align}
&(a^*b^*ba)\xi(x,y,r)  \notag \\
&=\int\overline{a(-\tilde{x},-\tilde{y},r)}\bar{e}\bigl[\eta_{\lambda}(r)
\beta(\tilde{x},y)\bigr](e^{\lambda\tilde{r}})^n\overline{b\bigl(
e^{\lambda\tilde{r}}(x-\tilde{x}),e^{\lambda\tilde{r}}(y-\tilde{y}),
-\tilde{r}\bigr)}  \notag \\
&\qquad(e^{\lambda\hat{r}})^n b\bigl(e^{\lambda\tilde{r}}(x-\tilde{x}),
e^{\lambda\tilde{r}}(y-\tilde{y}),\hat{r}\bigr)a(\hat{x},\hat{y},
r-\tilde{r}-\hat{r})  \notag \\
&\qquad\bar{e}\bigl[\eta_{\lambda}(r-\tilde{r}-\hat{r})\beta(\hat{x},
e^{\lambda(\tilde{r}+\hat{r})}(y-\tilde{y})-\hat{y})\bigr]  \notag \\
&\qquad\xi\bigl(e^{\lambda(\tilde{r}+\hat{r})}(x-\tilde{x})-\hat{x},
e^{\lambda(\tilde{r}+\hat{r})}(y-\tilde{y})-\hat{y},r-\tilde{r}-\hat{r}
\bigr)\,d\tilde{x}d\tilde{y}d\tilde{r}d\hat{r}d\hat{x}d\hat{y}.
\notag
\end{align}
So if we let $(\xi_i)$ be an orthonormal basis in ${\mathcal H}$,
we would have, by using change of variables:
\begin{align}
&\nu(a^*b^*ba)=\operatorname{Tr}(a^*b^*ba)=\sum_i\bigl\langle(a^*b^*ba)
\xi_i,\xi_i\bigr\rangle  \notag \\
&=\int\overline{a(-\tilde{x},-\tilde{y},r)}\bar{e}\bigl[\eta_{\lambda}(r)
\beta(\tilde{x},y)\bigr]\overline{b\bigl(e^{\lambda\tilde{r}}(x-\tilde{x}),
e^{\lambda\tilde{r}}(y-\tilde{y}),-\tilde{r}\bigr)}  \notag \\
&\qquad b\bigl(e^{\lambda\tilde{r}}(x-\tilde{x}),e^{\lambda\tilde{r}}
(y-\tilde{y}),-\tilde{r}\bigr)a(-\tilde{x},-\tilde{y},r)\bar{e}\bigl[
\eta_{\lambda}(r)\beta(-\tilde{x},y)\bigr]\,d\tilde{x}d\tilde{y}d\tilde{r}
dxdydr  \notag  \\
&=\int\overline{a(\tilde{x},\tilde{y},r)}a(\tilde{x},\tilde{y},r)
(e^{-2\lambda\tilde{r}})^n\overline{b(x,y,\tilde{r})}b(x,y,\tilde{r})\,
d\tilde{x}d\tilde{y}d\tilde{r}dxdydr   \notag \\
&=\varphi(a^*\times_A a)\hat{\varphi}(b^*\times_{\hat{A}}b).  \notag
\end{align}

Since ${\mathcal A}$ and $\hat{\mathcal A}$ generate the von Neumann
algebras $M_A$ and $M_{\hat{A}}$, while $M_A M_{\hat{A}}$ is dense in
${\mathcal B}({\mathcal H})$, this characterizes $\nu$.
\end{proof}

 The implication of this proposition is that for some well-chosen
element $a\in A$, the map $b\mapsto\nu(a^*ba)$ is a scalar multiple
of the weight $\hat{\varphi}(b)$. So proving the left invariance of
$\hat{\varphi}$ will be equivalent to showing the left invariance of
$\nu(a^*\,\cdot\,a)$.  This is done in Theorem \ref{li} below, with
a short lemma preceding it.  The steps are very similar to the proof
of Theorem 3.9 of \cite{BJKppha}.

\begin{lem}\label{lemmawk}
Let $(\xi_l)$ be an orthonormal basis for ${\mathcal H}$.  For $\zeta
\in{\mathcal H}$, consider the element $w_k=(\omega_{\zeta,\xi_k}
\otimes\operatorname{id})(U_A)\in{\mathcal B}({\mathcal H})$.  Then
we have:
$$
\sum_k\langle w_k\xi_l,w_k\xi_j\rangle=\langle\zeta,\zeta\rangle
\langle\xi_l,\xi_j\rangle.
$$
\end{lem}

\begin{rem}
The well-known definition of the forms of the type $\omega_{\zeta,\xi}
\in{\mathcal B}({\mathcal H})_*$ was given in the proof of Proposition
\ref{rho}.  Meanwhile, we know from \eqref{(Lomega)} that $w_k\in A$
(See also \cite{BS} and Proposition 2.2 of \cite{BJKppha}.).
\end{rem}

\begin{proof}
We take advantage of the fact that $(\xi_k)$ is an orthonormal basis,
and that $U_A$ is a unitary operator.  We have:
$$
\sum_k\langle w_k\xi_l,w_k\xi_j\rangle=\bigl\langle U_A(\zeta\otimes\xi_l),
U_A(\zeta\otimes\xi_j)\bigr\rangle=\langle\zeta\otimes\xi_l,\zeta\otimes
\xi_j\rangle=\langle\zeta,\zeta\rangle\langle\xi_l,\xi_j\rangle.
$$
\end{proof}

\begin{theorem}\label{li}
For any positive element $b\in\hat{A}$ such that $\hat{\varphi}(b)
<\infty$, and for positive $\omega\in{\hat{A}}^*$, we have:
$$
\hat{\varphi}\bigl((\omega\otimes\operatorname{id})(\hat{\Delta}b)\bigr)
=\omega(1)\hat{\varphi}(b).
$$
\end{theorem}

\begin{proof}
As suggested above in our comments following Proposition \ref{nutrace},
we may prove this for $\nu(a^*\,\cdot\,a)$, where $a\in{\mathcal A}
(\subseteq{\mathfrak N}_{\varphi})$ is a fixed element.

Let $b\in{{\mathfrak M}_{\hat{\varphi}}}^+$ and let $\omega\in{A^*}_+$.
Without loss of generality, we may assume that $\omega$ is the vector
state of the form $\omega=\omega_{\zeta,\zeta}$, for $\zeta\in{\mathcal H}$
We then have:
\begin{align}
(\omega\otimes\operatorname{id})(\hat{\Delta}b)&=(\omega_{\zeta,\zeta}
\otimes\operatorname{id})\bigl({U_A}^*(1\otimes b)U_A\bigr)  \notag \\
&=\sum_k\bigl[(\omega_{\xi_k,\zeta}\otimes\operatorname{id})({U_A}^*)\bigr]
b\bigl[(\omega_{\zeta,\xi_k}\otimes\operatorname{id})(U_A)\bigr]
=\sum_k w_k^*b^{\frac{1}{2}}b^{\frac{1}{2}}w_k.   \notag
\end{align}
Here $(\xi_k)$ is an orthonormal basis for ${\mathcal H}$, and the sums above
are convergent in the $\sigma$-weak topology on $M_{\hat{A}}$ (See Lemma 3.8
of \cite{BJKppha}.).  Also for convenience, we wrote $w_k=(\omega_{\zeta,\xi_k}
\otimes\operatorname{id})(U_A)$.

Let us use the result of the previous lemma and calculate:
\begin{align}
&\nu\bigl(a^*(\omega_{\zeta,\zeta}\otimes\operatorname{id})(\hat{\Delta}b)a
\bigr)=\sum_k\nu(a^*w_k^*b^{\frac{1}{2}}b^{\frac{1}{2}}w_ka)
\notag \\
&=\sum_{k,l}\langle w_kb^{\frac{1}{2}}a\xi_l,w_kb^{\frac{1}{2}}a\xi_l\rangle
\qquad\quad {\text { $\nu:=\operatorname{Tr}$ [trace on ${\mathcal B}
({\mathcal H})$]}}  \notag \\
&=\sum_l\langle\zeta,\zeta\rangle\langle b^{\frac{1}{2}}a\xi_l,b^{\frac{1}{2}}a
\xi_l\rangle\qquad\quad {\text { by Lemma \ref{lemmawk}}}  \notag \\
&=\langle\zeta,\zeta\rangle\operatorname{Tr}(a^*ba)
=\|\omega\|\nu(a^*ba)=\omega(1)\nu(a^*ba).    \notag
\end{align}
Since $\nu(a^*ba)=\varphi(a^*a)\hat{\varphi}(b)$, and since $\varphi(a^*a)$
is a positive constant, this will give us the proof that $\hat{\varphi}$ is
left invariant.
\end{proof}

 The left invariance we have just verified is a weak form, but by general
theory \cite{KuVa}, this is actually sufficient.  This establishes the
proof that $\hat{\varphi}$ is a legitimate Haar weight for $(\hat{A},\hat
{\Delta})$, in the sense that it is a proper, faithful, KMS weight which
is left invariant.  In our case, unlike the case of ${\varphi}_A$, the
Haar weight $\hat{\varphi}_{\hat{A}}$ is actually unimodular (Recall
Proposition \ref{Haarfunctional}).  Since this is the case, no extra
discussion is necessary for the ``right Haar weight'' or the ``modular
function''.  Summarizing the results of this section, we have the following
theorem:

\begin{theorem}
The pair $(\hat{A},\hat{\Delta})$, together with its additional
structures including the antipode and the (unimodular) Haar weight,
is a $C^*$-algebraic locally compact quantum group, in the sense
of Kustermans and Vaes.
\end{theorem}

 As we have made our case throughout this section , we may now regard
$(\hat{A},\hat{\Delta})$ as the {\em quantum Heisenberg group\/} (i.\,e.
``quantized $C_0(H)$'').  On the other hand, we remark here that our
example is different (and slightly more complicated) from the earlier
example of a quantum Heisenberg group obtained by Van Daele \cite{VD}.

 For instance, the Haar weight in the earlier example (although it was
not explicitly constructed in that paper) is a trace, while ours is
non-tracial.  The dual object of Van Daele's example is the example
by Rieffel \cite{Rf5}, while in our case, $(A,\Delta)$ of \cite{BJKp2},
\cite{BJKppha} plays that role.  These differences can be understood
more clearly if we consider the classical limits and compare the Poisson
structures: Our examples $(A,\Delta)$ and $(\hat{A},\hat{\Delta})$
were obtained by quantizing a certain non-linear Poisson structure,
while the examples of Rieffel's (\cite{Rf5}) and Van Daele's (\cite{VD})
correspond to a linear Poisson structure.

\section{Duality}

 The relationship between our two quantum groups $(A,\Delta)$ and $(\hat{A},
\hat{\Delta})$ is essentially the same as the relationship between $C^*(H)$
and $C_0(H)$.  Actually, it is a general fact that given a locally compact
quantum group $(B,\Delta)$, one can construct the dual quantum group $(\hat{B},
\hat{\Delta})$ within the category of locally compact quantum groups, and
that the generalized Pontryagin-type duality holds: That is, $(\hat{\hat{B}},
\hat{\hat{\Delta}})\cong (B,\Delta)$.  Refer \cite{MN}, \cite{KuVa} for the
general discussion on the duality of locally compact quantum groups.

 Our goal in this section is to see how the general theory is reflected in
the case of our specific examples.  Most of the results below are more or
less obvious and are direct consequences of general theory.  On the other
hand, several of these will be useful in our future applications.

\subsection{The dual pairing between ${\mathcal A}$ and $\hat{\mathcal A}$.}
Our quantum groups $(A,\Delta)$ and $(\hat{A},\hat{\Delta})$ are obtained
as two Hopf $C^*$-algebras associated with the multiplicative unitary operator
$U_A$ (as in \cite{BS}).  But unlike in the case of (finite-dimensional) Hopf
algebras, we do not actually have a dual pairing at the level of $C^*$-algebras
$A$ and $\hat{A}$.  What we do have is the dual pairing at the dense function
algebra level of ${\mathcal A}$ and $\hat{\mathcal A}$.  This is described
in the following proposition.

\begin{prop}
\begin{enumerate}
\item The dual pairing exists between ${\mathcal A}$ and $\hat{\mathcal A}$
such that for $f(=L_f)\in{\mathcal A}$ and $g(=\rho_g)\in\hat{\mathcal A}$,
we have:
$$
\langle f,g\rangle=\int f(x,y,r)g(e^{\lambda r}x,e^{\lambda r}y,-r)\,dxdydr.
$$
This is equivalent to the following pairing suggested by the multiplicative
unitary operator:
$$
\bigl\langle L(\omega),\rho({\omega}')\bigr\rangle=(\omega\otimes{\omega}')
(U_A)=\omega\bigl(\rho({\omega}')\bigr)={\omega}'\bigl(L(\omega)\bigr),
$$
where $L(\omega)=(\omega\otimes\operatorname{id})(U_A)\in A$ and $\rho
({\omega}')=(\operatorname{id}\otimes{\omega}')(U_A)\in\hat{A}$ are as
defined earlier with $\omega,{\omega}'\in{\mathcal B}({\mathcal H})_*$.
\item The dual pairing given above is compatible with the Hopf algebra
structures on ${\mathcal A}$ and $\hat{\mathcal A}$.  Indeed, for $f,
f_1,f_2\in{\mathcal A}$ and $g,g_1,g_2\in\hat{\mathcal A}$ (so $f=L_f$,
$g=\rho_g$, ...), we have:
\begin{align}
&\bigl\langle f_1\otimes f_2,\hat{\Delta}(g)\bigr\rangle = \bigl\langle
f_1\times_{\mathcal A}f_2,g\bigr\rangle,  \notag \\
&\bigl\langle f,g_1\times_{\hat{\mathcal A}}g_2\bigr\rangle=\bigl\langle
\Delta(f),g_1\otimes g_2\bigr\rangle,  \notag \\
&\bigl\langle S(f),g\bigr\rangle=\bigl\langle f,\hat{S}(g)\bigr\rangle,
\quad\bigl\langle f,g^*\bigr\rangle=\overline{\bigl\langle S(f)^*,
g\bigr\rangle}.  \notag
\end{align}
\end{enumerate}
\end{prop}

\begin{proof}
As described in the first part of the proposition, our definition of
the dual pairing was suggested by the theory of multiplicative unitary
operators.  For this, we use the same technique as in the proofs of
Propositions \ref{rho} and \ref{antipode}.  That is, consider $\omega_{\xi,
\eta}$, with $\xi,\eta$ being continuous functions with compact support
(contained in ${\mathcal H}$), so that we can realize the expressions like
$(\operatorname{id}\otimes\omega_{\xi,\eta})(U_A)$ as continuous functions
having compact support.

Once we take the above definition as our dual pairing, the verification
of the statements in the second part is very much straightforward.  All
we need to do is to remember the expressions of various operations
(for instance, equations \eqref{(Aops)}, \eqref{(Acoproduct)}, \eqref
{(multiplication)}, \eqref {(involution)} and Propositions \ref
{comultiplication} and \ref{antipode}) and just carry out the calculations.
For the first relation:
\begin{align}
&\bigl\langle f_1\otimes f_2,\hat{\Delta}(g)\bigr\rangle  \notag \\
&=\int f_1(x,y,\tilde{r})f_2(x',y',\tilde{r})g(e^{\lambda\tilde{r}}x
+e^{\lambda\tilde{r}}x',e^{\lambda\tilde{r}}y +e^{\lambda\tilde{r}}y',
-\tilde{r})  \notag \\
&\qquad e\bigl[\eta_{\lambda}(-\tilde{r})\beta(e^{\lambda\tilde{r}}x,
e^{\lambda\tilde{r}}y')\bigr]\,dxdydx'dy'd\tilde{r}  \notag \\
&=\int f_1(x,y,\tilde{r})f_2(x'-x,y'-y,\tilde{r})\bar{e}\bigl[\eta_{\lambda}
(\tilde{r})\beta(x,y'-y)\bigr]  \notag \\
&\qquad g(e^{\lambda\tilde{r}}x',e^{\lambda\tilde{r}}y',
-\tilde{r})\,dxdydx'dy'd\tilde{r}  \notag \\
&=\int(f_1\times_A f_2)(x',y',\tilde{r})g(e^{\lambda\tilde{r}}x',e^{\lambda
\tilde{r}}y',-\tilde{r})\,dx'dy'd\tilde{r}=\langle f_1\times_A f_2,g\rangle.
\notag
\end{align}
The other relations can be verified similarly.  Note that except the one
involving the ${}^*$ operation, the relations are exactly the ones we see
from ordinary Hopf algebra theory.
\end{proof}

\subsection{Duality at the Poisson--Lie group level.}
We have not been much emphasizing the role of the Poisson geometry in this
paper, but a brief discussion of the classical limit counterparts would
be useful here.  We have been arguing that $(A,\Delta)$ is a ``quantized
$C^*(H)$'' (See remark at the end of section 1, as well as our previous
papers \cite{BJKhj}, \cite{BJKppdress}.).  And we saw throughout section 2
that it is reasonable to consider $(\hat{A},\hat{\Delta})$ as a ``quantized
$C_0(H)$''.

 Meanwhile, in \cite{BJKp2}, we have made our case that $(A,\Delta)$ is
also a ``quantized $C_0(G)$'', where $G$ is the dual Poisson--Lie group
of $H$.  For the case of $(\hat{A},\hat{\Delta})$, we can actually regard
it as a ``quantized $C^*(G)$''.  To illustrate just one aspect of this,
recall the formula for the product on $\hat{\mathcal A}$ as given in \eqref
{(multiplication)}.  If we express this at the level of functions in the
$(p,q,r)$ variables (again by using the partial Fourier transform), it
becomes:
$$
(f\times_{\hat{A}}g)(p,q,r)=\int(e^{-2\lambda\tilde{r}})^n f(\tilde{p},
\tilde{q},\tilde{r})g(e^{-\lambda\tilde{r}}p-e^{-\lambda\tilde{r}}\tilde{p},
e^{-\lambda\tilde{r}}q-e^{-\lambda\tilde{r}}\tilde{q},r-\tilde{r})\,d\tilde{p}
d\tilde{q}d\tilde{r}.
$$
But if we assume that $H^*$ has the group structure given by the
multiplication law:
$$
(p,q,r)(p',q',r')=(e^{\lambda r'}p+p',e^{\lambda r'}q+q',r+r'),
$$
which is exactly the multiplication law for the dual Poisson--Lie group
$G$ of $H$ as defined in \cite{BJKp2}, then the above expression for the
product on $\hat{\mathcal A}$ can be written as:
$$
(f\times_{\hat{A}}g)(p,q,r)=\int(e^{-2\lambda\tilde{r}})^n f(\tilde{p},
\tilde{q},\tilde{r})g\bigl((p,q,r)(\tilde{p},\tilde{q},\tilde{r})^{-1}
\bigr)\,d\tilde{p}d\tilde{q}d\tilde{r}.
$$
Since $(e^{-2\lambda\tilde{r}})^n\,d\tilde{p}d\tilde{q}d\tilde{r}$ is
the right Haar measure for the group $G(=H^*)$, this means that it
is really the convolution product.  In other words, we notice that
$\hat{A}\cong C^*(G)$ as a $C^*$-algebra, where $C^*(G)$ is realized
as an operator algebra via the right regular representation of $G$.

 These observations illustrate that the duality between $(A,\Delta)$
and $(\hat{A},\hat{\Delta})$ is the quantum counterpart to the
Poisson--Lie group duality between $G$ and $H$.  This point of view
is certainly very useful in any applications involving our quantum
groups.  The duality picture will be enhanced when we consider the
``quantum double'' of our examples (Just as the ``double Poisson--Lie
group'' $H\Join G$ and the ``dressing orbits'' play a useful role
\cite{BJKhj}, \cite{BJKppdress}.).  In a future paper, we will give
a discussion on the quantum double construction, again within the
framework of $C^*$-algebraic, locally compact quantum groups.

\subsection{The ``opposite'' and ``co-opposite'' Hopf $C^*$-algebras.}

 By slightly modifying our fundamental multiplicative unitary operator
$U_A$, we are able to construct a few different forms of the quantum
Heisenberg group and the quantum Heisenberg group algebra.  Borrowing
terminologies from Hopf algebra theory, they will more or less correspond
to ``opposite'' or ``co-opposite'' algebras, and ``opposite dual'' or
``co-opposite dual'' algebras.

 Let $j\in{\mathcal B}({\mathcal H})$ be defined by
$$
j\xi(x,y,r)=(e^{\lambda r})^n\bar{e}\bigl[\eta_{\lambda}(r)\beta(x,y)\bigr]
\xi(-e^{\lambda r}x,-e^{\lambda r}y,-r).
$$
Then $j$ is a unitary operator such that $j^2=1$.  Note that the operator
$j$ can be written as $j=\hat{J}J=J\hat{J}$, where $J$ and $\hat{J}$ are
the operators we saw earlier in our discussions on the antipode and the
${}^*$-operation.  Incorporating the operator $j$ to our fundamental
multiplicative unitary operator $U_A$, we obtain the following:

\begin{prop}\label{unitarymultiplicative}
The following operators are all regular multiplicative unitary operators
(in the sense of Baaj and Skandalis) in ${\mathcal B}({\mathcal H}\otimes
{\mathcal H})$.  Here $\Sigma$ denotes the flip.
\begin{alignat}{2}
U_A&  &\widehat{U_A}&=\Sigma(j\otimes1)U_A(j\otimes1)\Sigma  \notag \\
\widetilde{U_A}&=(j\otimes1)(\Sigma U_A \Sigma)(j\otimes1)\qquad\quad&
\widehat{\widehat{U_A}}&=\widetilde{\widetilde{U_A}}=(j\otimes j)U_A
(j\otimes j) \notag
\end{alignat}
The same is true of the operators of the form $\Sigma X^*\Sigma$ for
any of the multiplicative unitary operators $X$ above.
\end{prop}

\begin{rem}
The verification is a straightforward calculation.  What is really
going on is that the triple $({\mathcal H},U_A,j)$ forms a {\em Kac
system\/}, in the terminology of Baaj and Skandalis (See section 6
of \cite{BS}.).
\end{rem}

 We will obtain several $C^*$-bialgebras from these operators.  Before
we give descriptions of them, let us give the following definitions on
``opposite'' and ``co-opposite'' algebra/coalgebra structures.

\begin{defn}
\begin{enumerate}
\item On ${\mathcal A}$, define instead the ``opposite multiplication''
by $(f,g)\mapsto g\times_A f$.  We keep the same involution.  We will
denote this opposite algebra by ${\mathcal A}^{\operatorname{op}}$.
Similarly, we can define $\hat{\mathcal A}^{\operatorname{op}}$,
whose multiplication is given by $(f,g)\mapsto g\times_{\hat{A}}f$.
\item Define $R:{\mathcal A}^{\operatorname{op}}\to{\mathcal B}
({\mathcal H})$ by
$$
(R_f\xi)(x,y,r):=\int f(\tilde{x},\tilde{y},r)\xi(x-\tilde{x},
y-\tilde{y},r)\bar{e}\bigl[\eta_{\lambda}(r)\beta(x-\tilde{x},
\tilde{y})\bigr]\,d\tilde{x}d\tilde{y}.
$$
It is a ${}^*$-representation of ${\mathcal A}^{\operatorname{op}}$.
Actually, ${\mathcal A}^{\operatorname{op}}$ is a pre-$C^*$-algebra,
together with the $C^*$-norm $\|f\|:=\|R_f\|$.  We will denote by
$A^{\operatorname{op}}$ the $C^*$-algebra completion in ${\mathcal B}
({\mathcal H})$ of ${\mathcal A}^ {\operatorname{op}}$.  That is,
$A^{\operatorname{op}}=\overline {R({\mathcal A}^{\operatorname{op}})}^
{\|\ \|}$.
\item Define $\lambda:\hat{\mathcal A}^{\operatorname{op}}\to{\mathcal B}
({\mathcal H})$ by
$$
(\lambda_f\zeta)(x,y,r):=\int f(e^{\lambda\tilde{r}}x,e^{\lambda
\tilde{r}}y,r-\tilde{r})\zeta(x,y,\tilde{r})\,d\tilde{r}.
$$
It is a ${}^*$-representation of $\hat{\mathcal A}^{\operatorname{op}}$.
As above, we can define the $C^*$-algebra $\hat{A}^{\operatorname{op}}$
as $\hat{A}^{\operatorname{op}}=\overline{\lambda(\hat{\mathcal A}^
{\operatorname{op}})}^{\|\ \|}\bigl(\subseteq{\mathcal B}({\mathcal H})
\bigr)$.
\end{enumerate}
\end{defn}

\begin{rem}
The above definitions resemble the characterizations of the $C^*$-algebras
$A$ and $\hat{A}$ in ${\mathcal B}({\mathcal H})$.  And the roles played
by $R({\mathcal A}^{\operatorname{op}})$ and $\lambda(\hat{\mathcal A}^
{\operatorname{op}})$ are exactly the same ones played by $L({\mathcal A})$
and $\rho(\hat{\mathcal A})$.  Meanwhile, on a related note concerning
the enveloping von Neumann algebras, we have: $M_{{\mathcal A}^
{\operatorname{op}}}={M_A}'$, and $M_{\hat{\mathcal A}^{\operatorname{op}}}
={M_{\hat{A}}}'$.
\end{rem}

\begin{defn}
\begin{enumerate}
\item For the function $f$ contained in ${\mathcal A}$ (or in ${\mathcal A}^
{\operatorname{op}}$), define ${\Delta}^{\operatorname{cop}}f$ by
\begin{align}
&{\Delta}^{\operatorname{cop}}f(x,y,r;x',y',r')  \notag \\
&=\int f(x,y,r+r')\bar{e}\bigl[\tilde{p}\cdot(e^{\lambda r}x-x')
+\tilde{q}\cdot(e^{\lambda r}y-y')\bigr]\,d\tilde{p}d\tilde{q},  \notag
\end{align}
which is a Schwartz function having compact support in the $r$ and the $r'$
variables.
\item For $f$ contained in $\hat{\mathcal A}$ (or in $\hat{\mathcal A}
^{\operatorname{op}}$), let $\hat{\Delta}^{\operatorname{cop}}f$ be the
Schwartz function having compact support in $r$ and $r'$, defined by
\begin{align}
&\hat{\Delta}^{\operatorname{cop}}f(x,y,r;x',y',r')  \notag \\
&=\int
f(x+x',y+y',\tilde{r})e\bigl[\eta_{\lambda}(\tilde{r})\beta(x',y)\bigr]
e\bigl[\tilde{r}(z+z')\bigr]\bar{e}[zr+z'r']\,d\tilde{r}dzdz'.  \notag
\end{align}
\end{enumerate}
\end{defn}

\begin{rem}
As the names suggest, these are the ``co-opposite comultiplications'' 
(Compare the above definitions with our earlier definitions of $\Delta f$
and $\hat{\Delta}f$ given in \eqref{(Acoproduct)} and Proposition \ref
{comultiplication}.).  Indeed, for $f\in{\mathcal A}$, we would have:
$(L\otimes L)_{{\Delta}^{\operatorname{cop}}f}=(\chi\circ{\Delta})(L_f)$,
where $\chi$ is the flip.  Similar comment holds for $\hat{\Delta}^
{\operatorname{cop}}$.  Meanwhile, just as were the cases of $\Delta$ and
$\hat{\Delta}$, the above maps ${\Delta}^{\operatorname{cop}}$ and $\hat
{\Delta}^{\operatorname{cop}}$ can be also extended to the $C^*$-algebra
level (See Proposition \ref{opcopalgebras} below.).
\end{rem}

 Let us turn our attention back to the multiplicative unitary operators
in Proposition \ref{unitarymultiplicative}.  For each of the multiplicative
unitary operators $V$, we can consider $\overline{\bigl\{(\omega\otimes
\operatorname{id})(V):\omega\in{\mathcal B}({\mathcal H})_*\bigr\}}^{\|\ \|}$
(the ``left slices'') and $\overline{\bigl\{(\operatorname{id}\otimes\omega)
(V):\omega\in{\mathcal B} ({\mathcal H})_*\bigr\}}^{\|\ \|}$ (the ``right
slices'') contained in ${\mathcal B}({\mathcal H})$. They are described
below:

\begin{prop}\label{opcopalgebras}
For $f\in{\mathcal A}$ and $g\in\hat{\mathcal A}$, we have:
\begin{align}
&U_A(L_f\otimes1){U_A}^*={\widehat{U_A}}^*(1\otimes L_f)\widehat{U_A}
=(L\otimes L)_{{\Delta}f},  \notag \\
&\widehat{\widehat{U_A}}(R_f\otimes1){\widehat{\widehat{U_A}}}^*
={\tilde{U_A}}^*(1\otimes R_f)\tilde{U_A}
=(R\otimes R)_{{\Delta}^{\operatorname{cop}}f},  \notag \\
&{U_A}^*(1\otimes\rho_g)U_A=\tilde{U_A}(\rho_g\otimes1){\tilde{U_A}}^*
=(\rho\otimes\rho)_{\hat{\Delta}g}, \notag \\
&{\widehat{\widehat{U_A}}}^*(1\otimes\lambda_g)\widehat{\widehat{U_A}}
=\widehat{U_A}(\lambda_g\otimes1){\widehat{U_A}}^*=(\lambda\otimes
\lambda)_{\hat{\Delta}^{\operatorname{cop}}g}.  \notag
\end{align}
From this, we obtain the following results (Here, the comultiplications
are understood as defined at (extended to) the $C^*$-algebra level.):
\begin{enumerate}
\item $U_A$ determines two Hopf $C^*$-algebras $(A,\Delta)$ and $(\hat{A},
\hat{\Delta})$.  And $\Sigma{U_A}^*\Sigma$ determines $(\hat{A},\hat{\Delta}^
{\operatorname{cop}})$ and $(A,{\Delta}^{\operatorname{cop}})$.
\item $\widehat{U_A}$ determines $(\hat{A}^{\operatorname{op}},\hat{\Delta}^
{\operatorname{cop}})$ and $(A,\Delta)$, while $\Sigma{\widehat{U_A}}^*\Sigma$
determines $(A,{\Delta}^{\operatorname{cop}})$ and $(\hat{A}^{\operatorname
{op}},\hat{\Delta})$.
\item $\tilde{U_A}$ determines $(\hat{A},\hat{\Delta})$ and $(A^{\operatorname
{op}},\Delta^{\operatorname{cop}})$, while $\Sigma{\tilde{U_A}}^*\Sigma$
determines $(A^{\operatorname{op}},{\Delta})$ and $(\hat{A},\hat{\Delta}^
{\operatorname{cop}})$.
\item $\widehat{\widehat{U_A}}$ determines $(A^{\operatorname{op}},{\Delta}^
{\operatorname{cop}})$ and $(\hat{A}^{\operatorname{op}},\hat{\Delta}^
{\operatorname{cop}})$, while $\Sigma{\widehat{\widehat{U_A}}}^*\Sigma$
determines $(\hat{A}^{\operatorname{op}},\hat{\Delta})$ and $(A^{\operatorname
{op}},{\Delta})$.
\end{enumerate}
\end{prop}

\begin{proof}
Case (1) repeats the results of the sections 1 and 2.  And case (2) was
considered in Appendix (section 6) of \cite{BJKppha}.

To obtain the $C^*$-algebras corresponding to the multiplicative unitary
operators, we adopt the method we used in Proposition \ref{rho}.  Their
comultiplications can be read from one of the equations in the first part,
which can be proved by a straightforward calculation (See also Proposition
6.8 of \cite{BS}.).  In addition, these equations justify our viewing the
comultiplications as defined at the $C^*$-algebra level.
\end{proof}

\begin{cor}
Let the notation be as in Proposition \ref{opcopalgebras}.  We have:
\begin{enumerate}
\item $U_A\in M(\hat{A}\otimes A)$.
\item $\widehat{U_A}\in M(A\otimes\hat{A}^{\operatorname{op}})$.
\item $\tilde{U_A}\in M(A^{\operatorname{op}}\otimes\hat{A})$.
\item $\widehat{\widehat{U_A}}\in M(\hat{A}^{\operatorname{op}}
\otimes A^{\operatorname{op}})$.
\end{enumerate}
\end{cor}

\begin{proof}
The results are immediate consequences of Proposition \ref{opcopalgebras}.
Refer to Proposition 3.6 of \cite{BS}.
\end{proof}

 At this moment, what we have in Proposition \ref{opcopalgebras} are just
$C^*$-bialgebras.  But together with the appropriate Haar weights, they become
locally compact quantum groups (We may work with the same Haar weights that
we have used for $(A,\Delta)$ and $(\hat{A},\hat{\Delta})$.  In some cases,
the roles of left Haar weight and the right Haar weight have to be reversed.).
In this way, we obtain various versions of the quantum Heisenberg group
and its dual (On the other hand, note that we have $(A,\Delta)\cong
(A^{\operatorname{op}}, {\Delta}^{\operatorname{cop}})$, via the antipode.).

%%%%%%%%%%%%%%%%%%%%%%%%%%%%%%%%%%%%%%%%%%%%%%%%%%%%%%%%%%%%%%%%%%%%%%%%%

%\pagebreak

\bibliography{ref}

\providecommand{\bysame}{\leavevmode\hbox to3em{\hrulefill}\thinspace}
\providecommand{\MR}{\relax\ifhmode\unskip\space\fi MR }
% \MRhref is called by the amsart/book/proc definition of \MR.
\providecommand{\MRhref}[2]{%
  \href{http://www.ams.org/mathscinet-getitem?mr=#1}{#2}
}
\providecommand{\href}[2]{#2}
\begin{thebibliography}{10}

\bibitem{BS}
S.~Baaj and G.~Skandalis, \emph{Unitaires multiplicatifs et dualit\'e pour les
  produits crois\'es de {$C^*$}-alg\`ebres}, Ann. Scient. \'Ec. Norm. Sup.,
  $4^e$ s\'erie \textbf{t. 26} (1993), 425--488 (French).

\bibitem{Cm1}
F.~Combes, \emph{Poids sur une {$C^*$}-alg\`ebre}, J. Math. Pures et Appl.
  \textbf{47} (1968), 57--100 (French).

\bibitem{Cm2}
\bysame, \emph{Poids associ\'e \`a une alg\`ebre hilbertienne \`a gauche},
  Compos. Math. \textbf{23} (1971), 49--77 (French).

\bibitem{ES}
M.~Enock and J.~M. Schwartz, \emph{Kac {A}lgebras and {D}uality of {L}ocally
  {C}ompact {G}roups}, Springer-Verlag, 1992.

\bibitem{BJKp2}
B.~J. Kahng, \emph{Non-compact quantum groups arising from {H}eisenberg type
  {L}ie bialgebras}, J. Operator Theory \textbf{44} (2000), 303--334.

\bibitem{BJKhj}
\bysame, \emph{${}^*$-representations of a quantum {H}eisenberg group algebra},
  Houston J. Math. \textbf{28} (2002), 529--552.

\bibitem{BJKppdress}
\bysame, \emph{Dressing orbits and a quantum {H}eisenberg group algebra}, 2002,
  preprint (accepted to appear in Illinois J. Math., available as
  math.OA/0211003 at http://lanl.arXiv.org).

\bibitem{BJKppha}
\bysame, \emph{Haar measure on a locally compact quantum group}, J. Ramanujan
  Math. Soc. \textbf{18} (2003), 385--414.

\bibitem{Ku1}
J.~Kustermans, \emph{{KMS}-weights on {$C^*$}-algebras}, 1997, preprint, Odense
  Universitet (funct-an/9704008).

\bibitem{KuVa}
J.~Kustermans and S.~Vaes, \emph{Locally compact quantum groups}, Ann. Scient.
  \'Ec. Norm. Sup., $4^e$ s\'erie \textbf{t. 33} (2000), 837--934.

\bibitem{MN}
T.~Masuda and Y.~Nakagami, \emph{A von {N}eumann algebra framework for the
  duality of the quantum groups}, Publ. RIMS, Kyoto Univ. \textbf{30} (1994),
  no.~5, 799--850.

\bibitem{Rf3}
M.~A. Rieffel, \emph{Lie group convolution algebras as deformation
  quantizations of linear {P}oisson structures}, Amer. J. Math. \textbf{112}
  (1990), 657--685.

\bibitem{Rf5}
\bysame, \emph{Some solvable quantum groups}, Operator Algebras and Topology
  (W.~B. Arveson, A.~S. Mischenko, M.~Putinar, M.~A. Rieffel, and S.~Stratila,
  eds.), Proc. OATE2 Conf: Romania 1989, Pitman Research Notes Math., no. 270,
  Longman, 1992, pp.~146--159.

\bibitem{St}
S.~Stratila, \emph{Modular {T}heory in {O}perator {A}lgebras}, Abacus Press,
  1981.

\bibitem{SZ}
I.~Szymczak and S.~Zakrzewski, \emph{Quantum deformations of the {H}eisenberg
  group obtained by geometric quantization}, J. Geom. Phys. \textbf{7} (1990),
  553--569.

\bibitem{VV}
S.~Vaes and L.~Vainerman, \emph{Extensions of locally compact quantum groups
  and the bicrossed product construction}, Adv. Math. \textbf{175} (2003),
  1--101.

\bibitem{Va}
J.~Vallin, \emph{{$C^*$}-alg\`ebres de {H}opf et {$C^*$}-alg\`ebres de {K}ac},
  Proc. London Math. Soc. \textbf{50} (1985), 131--174 (French).

\bibitem{VD}
A.~{Van Daele}, \emph{Quantum deformation of the {H}eisenberg group},
  Proceedings of the Satellite Conference of ICM-90, World Scientific,
  Singapore, 1991, pp.~314--325.

\bibitem{VDaxb}
\bysame, \emph{The {H}aar measure on some locally compact quantum groups},
  2001, preprint (available as math.OA/0109004 at http://lanl.arXiv.org).

\bibitem{VDoamp}
\bysame, \emph{The {H}eisenberg commutation relations, commuting squares and
  the {H}aar measure on locally compact quantum groups}, 2001, preprint (to
  appear in Proceedings of the OAMP Conference, Constantza, 2001).

\bibitem{Wr7}
S.~L. Woronowicz, \emph{From multiplicative unitaries to quantum groups}, Int.
  J. Math. \textbf{7} (1996), no.~1, 127--149.

\end{thebibliography}

\bibliographystyle{amsplain}

\end{document}